\documentclass{amsart}
\usepackage{hyperref}

\date{\today}

\newcommand{\Z}{{\mathbb Z}}
\newcommand{\R}{{\mathbb R}}

\newcommand{\be}{\begin{equation}}
\newcommand{\ee}{\end{equation}}

\newcommand{\ti}{\tilde}

\newcommand{\E}{\mathrm{e}}

\newcommand{\tr}{\mathrm{tr}}

\newcommand{\eps}{\varepsilon}

\newtheorem{theorem}{Theorem} [section]

\newtheorem{lemma}[theorem]{Lemma}
\newtheorem{prop}[theorem]{Proposition}

\numberwithin{equation}{section}

\begin{document}

\title[Optimality of log H\"older continuity of the IDS]{Optimality of log H\"older continuity of the integrated density of states}

\author[Z.\ Gan]{Zheng Gan}

\address{Department of Mathematics, Rice University, Houston, TX~77005, USA}
\email{\href{mailto:zheng.gan@rice.edu}{zheng.gan@rice.edu}}
\urladdr{\href{http://math.rice.edu/~zg2/}{http://math.rice.edu/$\sim$zg2/}}

\author[H.\ Kr\"uger]{Helge Kr\"uger}

\address{Department of Mathematics, Rice University, Houston, TX~77005, USA}
\email{\href{mailto:helge.krueger@rice.edu}{helge.krueger@rice.edu}}
\urladdr{\href{http://math.rice.edu/~hk7/}{http://math.rice.edu/$\sim$hk7/}}

\thanks{H.\ K.\ was supported by NSF grant DMS--0800100.}

\date{\today}


\begin{abstract}
 We construct examples, that log H\"older continuity of the integrated density
 of states cannot be improved. Our examples are limit-periodic.
\end{abstract}

\maketitle

\section{Introduction}

We investigate optimality of the log H\"older continuity of the
integrated density of states. Let $(\Omega,\mu)$ be a probability space,
$T:\Omega\to\Omega$ an invertible ergodic transformation, and $f:
\Omega\to\R$ a bounded measurable function. Define a potential
$V_{\omega}(n) = f(T^n\omega)$. The Sch\"odinger
operator $H_{\omega} :\ell^2(\Z) \to \ell^2(\Z) $ is defined
by
\be \label{equ:oper}
 H_{\omega} u(n) = u(n+1) + u(n-1) +V_{\omega}(n) u(n),
\ee
and the integrated density of states $k$ by
\be\label{eq:defids}
 k(E) = \lim_{N \to \infty} \int_{\Omega} \left(\frac{1}{N} \tr(P_{(-\infty,E)}(H_{\omega, [0, N-1]}))\right) d\mu(\omega),
\ee
where $H_{\omega,[0,N-1]}$ denotes the restriction of
$H_{\omega}$ to $\ell^2([0,N-1])$. Using the \textit{Thouless
formula}, Craig and Simon showed that

\begin{theorem}[Craig and Simon, \cite{cs2}]
 There exists a constant $C = C(\|f\|_{\infty})$ such that
 \be
  |k(E) - k(\ti{E})| \leq \frac{C}{\log|E - \ti{E}|^{-1}}
 \ee
 for $|E - \ti{E}| \leq \frac{1}{2}$.
\end{theorem}

This is what is well known as \textit{log H\"older continuity}. We
will be interested in the optimality of this statement in the sense
that $\eps \mapsto \frac{1}{\log(\eps^{-1})}$ cannot be replaced by
another function, which goes to zero faster. It was shown by
Craig in \cite{cr}, that the regularity cannot be improved
to
$$
 \eps \mapsto \frac{1}{\log(\eps^{-1}) \log(\log(\eps^{-1}))^\beta},
$$
where $\beta > 1$. However, in the case of specific dynamical systems $(\Omega,\mu,T)$,
there exist many results, which improve the Craig--Simon result.
We just mention two. For quasi-periodic Schr\"odinger operators,
Goldstein and Schlag have shown in \cite{gs} that the integrated
density of states is H\"older continuous and computed the H\"older exponent,
and shown that the integrated density of states is almost everywhere
Lipschitz. For random Schr\"odinger operators, the integrated
density of states is even everywhere Lipschitz. This is known
as the Wegner estimate which can be found for example
in the exposition of Kirsch in \cite{dkkkr}.

Our interest in the question of optimality of the Craig--Simon results
comes from the importance of the Wegner estimate in multiscale analysis
(see the exposition of Kirsch). If one could improve the result to
a continuity of the form
$$
 \eps \mapsto \frac{1}{\log(\eps^{-1})^\beta}
$$
for some large enough $\beta > 1$, one would be able to use this for multiscale
analysis (see for example Theorem~3.12 in \cite{k}). Already Craig's
result shows that this is impossible, however one could hope that
a combination of an improved continuity result and an improvement of
multiscale analysis might remove the Wegner estimate assumption.
However, we will show that the continuity of integrated density
of states cannot be improved for all potentials beyond log H\"older
continuity.

A potential $V \in \ell^\infty(\Z)$ is called almost-periodic,
if the closure $\Omega$ of its translates is compact in the $\ell^\infty$ norm.
Furthermore, then $\Omega$ can be made into a compact group, with
an unique invariant Haar measure. For these our previous definition
of the integrated density of states \eqref{eq:defids} can be replaced
by
\be\label{eq:defids2}
 k_V(E) = \lim_{N \to \infty} \frac{1}{N} \tr(P_{(-\infty,E)}(H_{[0,N-1]})),
\ee
where $H_{[0,N-1]}$ denotes now the restriction of $\Delta + V$
to $\ell^2([0,N-1])$ with $\Delta$ the discrete Laplacian. This can
be found for example as Theorem~2.9 in Avron--Simon \cite{as}.

Next, $V$ is called $p$ periodic, if its $p$-th translate is equal to $V$.
Furthermore, $V$ is limit-periodic if it is the limit in the $\ell^\infty$
norm of periodic potentials. We denote by $\sigma(\Delta+V)$ the spectrum
of the operator $\Delta  + V$.

\begin{theorem}\label{thm:optimal2}
 Given any increasing continuous function $\varphi : \R^{+} \to \R^{+}$ with
 \be
  \lim_{x\to 0} \varphi(x) = 0
 \ee
 and a constant $C_0 > 0$, there is a limit-periodic $V$ satisfying
 $\|V\|_{\infty} \leq C_0$ such that its integrated density of states
 satisfies
 \be \label{equa:thm}
  \limsup_{E \to E_0} \frac{|k_V(E)  - k_V(E_0)| \log(|E - E_0|^{-1})}{\varphi(|E - E_0|)} =
  \infty,
 \ee
 for any $E_0 \in \sigma(\Delta  + V)$.
\end{theorem}

This result tells us, that with $\varphi$ as in the previous theorem, we cannot
have
$$
 |k_V(E)  - k_V(E_0)| \leq C \cdot \frac{\varphi(|E - E_0|)}{\log(|E-E_0|^{-1})}
$$
for any $C > 0$ and all $V$.
The proof of this theorem essentially happens in two parts. Given a periodic $V_0$ and $\eps > 0$
satisfying $\|V_0\| \leq C_0 - \eps$, we construct
a sequence $V_j$ of periodic potentials, with the following properties
\begin{enumerate}
 \item $V_j$ is $p_j$-periodic.
 \item The Lebesgue measure of $\sigma(\Delta + V_j)$
  \be
   \eps_j = |\sigma(\Delta + V_j)|
  \ee
  satisfies
  \be\label{eq:sizeepsj}
   \log(\eps_j^{-1}) \geq p_{j-1} \cdot p_{j} \cdot \varphi(2 \eps_j).
  \ee
 \item We have that
  \be\label{eq:Vjconv}
   \|V_j - V_{j-1}\| \leq \frac{\min(\eps, \eps_1, \dots, \eps_{j-1})}{2^j}.
  \ee
\end{enumerate}
Here $\| . \|$ denotes the $\ell^{\infty}$ norm.
The construction of these $V_j$ will be given in the next section
and uses the tools developed by Avila in \cite{a}.
Before proceeding with the proof of Theorem~\ref{thm:optimal2},
recall that positivity of the trace implies that
$$
 k_W(E - \|V - W\|) \leq k_V(E) \leq k_W(E + \|V-W\|)
$$
for any potentials $V$ and $W$.

\begin{proof}[Proof of Theorem~\ref{thm:optimal2}]
 By \eqref{eq:Vjconv}, we see that there exists a limiting potential $V$, such that
 for each $j$, we have that
 $$
  \|V_j - V\| \leq \eps_j = |\sigma(\Delta + V_j)|.
 $$
 Furthermore, since $\|V_0 - V\| \leq \eps$, we have that $\|V\|\leq C_0$.

 Next, fix $E_0 \in \sigma(\Delta + V)$ and let $j \geq 1$. By the previous equation, we have that
 there exists $E_1 \in \sigma(\Delta + V_j)$ such that
 $$
  |E_0 - E_1| \leq \eps_j = |\sigma(\Delta + V_j)|.
 $$
 Denote the band of $\sigma(\Delta + V_j)$ containing $E_1$ by $[E_-, E_+]$.
 By a general fact about periodic Schr\"odinger operators, we know that
 $$
  k_{V_j}(E_+) - k_{V_j}(E_-) = \frac{1}{p_j}.
 $$
 We thus get that
 $$
  k_V(E_+ + \eps_j) - k_V(E_- - \eps_j) \geq \frac{1}{p_j}
 $$
 Furthermore, the interval $[E_- - \eps_j, E_+ + \eps_j]$ contains
 $E_0$ and we can choose $E_j \in \{E_- - \eps_j, E_+ + \eps_j\}$ such that
 $$
  |k_V(E_0) - k_V(E_j)| \geq \frac{1}{2 p_j}
 $$
 and
 $$
  |E_0 - E_j| \leq 2 \eps_j.
 $$
 This implies the claim by \eqref{eq:sizeepsj}, since $j$ was arbitrary.
\end{proof}

One can slightly improve the above theorem, by for example showing that
there is not only one $V$, that satisfies the conclusion, but that in fact
the set is dense in the limit-periodic operators. However, we have not done
so, to keep the statement as simple as possible.

\section{Construction of the periodic potentials}

We will need the machinery developed by Avila in \cite{a}, in order
to prove our results. In the following, we let $\Omega$ be a totally
disconnected compact group, known as {\em Cantor group}. We furthermore
let $T: \Omega\to\Omega$ be a minimal translation on this group.
There is a decreasing sequence of Cantor subgroups
$$
 X_1 \supseteq X_2 \supseteq \dots
$$
such that the quotients
$$
 \Omega / X_k
$$
contain $p_k$ elements. We let $P_k$ be the subset of the continuous
functions $C(\Omega)$ on $\Omega$, which only depend on $\Omega/ X_k$.
$f$ is called $n$-periodic if $f(T^n \omega) = f(\omega)$ for every $\omega
\in \Omega$.  The elements of $P_k$ will be $p_k$ periodic.

We now fix $\omega \in \Omega$. We have that $\{f(T^n\omega)\}_{n\in\Z} \in \ell^{\infty}(\Z)$ is
limit-periodic, since the periodic $f$ are dense in $C(\Omega)$.
For a finite subset $F$ of the periodic potentials $P = \bigcup_{k \geq 1} P_k$,
we introduce the averaged Lyapunov exponent $L(E, F)$ as
\be
 L(E, F) = \frac{1}{\# F} \sum_{f \in F} L(E, f),
\ee
where $\# F$ denotes the number of elements of $F$ (with
multiplicities) and
$$
 L(E, f) = \lim_{N\to\infty} \frac{1}{N} \log\left\|\prod_{n=N}^{1} \begin{pmatrix} f(T^n \omega) - E & - 1 \\ 1 & 0 \end{pmatrix}\right\|
$$
is the Lyapunov exponent of the periodic potential. For $f\in C(\Omega)$,
we denote by $\Sigma(f)$ the spectrum of the operator $\Delta + f(T^n \omega)$.
We will use the following two lemmas of Avila \cite{a}, see also \cite{dg2}.

\begin{lemma}[Lemma~3.1. in \cite{a}]\label{lem:avila1}
 Let $B$ be an open ball in $C(\Omega)$,  let $F \subset P \cap B$
 be finite, and  let $0 < \eps < 1$. Then there exists a sequence $F_K
 \subset P \cap B$ such that
 \begin{enumerate}
  \item $L(E, F_K) > 0$ whenever $E \in \R$,
  \item $L(E, F_K) \rightarrow L(E, F)$ uniformly on compacts.
 \end{enumerate}
\end{lemma}

\begin{lemma}[Lemma~3.2. in \cite{a}]\label{lem:avila2}
 Let $B$ be an open ball in $C(\Omega)$, and let $F \subset P_k \cap
 B$ be a finite family of sampling functions. Then for every
 $N \ge 2$ and $K$ sufficiently large, there exists $F_K \subset P_K \cap B $ such that
 \begin{enumerate}
  \item $L(E, F_K) \rightarrow L(E, F)$ uniformly on compacts,
  \item The diameter of $F_K$ is at most $p_K^{-10}$,
  \item For every $\lambda \in \R$, if
    \be
     \inf_{E \in \R} L(E, F) \ge \delta \#F p_k,
    \ee
   then for every $f \in F_K$, the spectrum $\Sigma(f)$
   has Lebesgue measure at most $e^{-\delta p_K /2}$.
 \end{enumerate}
\end{lemma}


The construction of the $V_j$ will be accomplished by

\begin{prop} \label{prop:optimal}
 Given a continuous function $\psi: \R^+ \to \R^+$ satisfying
 \be
  \lim_{x \to 0} \psi(x) = 0,
 \ee
 and $p$-periodic $f$ and $\eps > 0$, then there
 exists a $\ti{p}$-periodic function $\ti{f}$, such that
 \be
  \|f - \ti{f}\|_{\infty} \leq \eps,
 \ee
 and
 \be \label{equa:opt}
  \log(|\Sigma(\ti{f})|^{-1}) \geq \ti{p} \cdot \psi(|\Sigma(\ti{f})|).
 \ee
\end{prop}

\begin{proof}
 By Lemma~\ref{lem:avila1}, we can
 find a finite family of $p_1$-periodic potentials $F_1$ within
 $B_{\frac{1}{2} \eps}(f)$ such that
 $$
  \delta_1 = \frac{1}{\# F_1 p_1}L(E, F_1) = \frac{1}{\# F_1 p_1} \cdot\frac{1}{\# F_1} \sum_{f \in F_1} L(E,f) > 0.
 $$
 Applying Lemma~\ref{lem:avila2} to $F_1$, we can get a finite family
 $F_2 \subset B_{\eps}(f)$ of $\tilde{p}$-periodic potentials, where
 we might require $\ti{p}$ to be arbitrarily large. Let $\ti{f}$ be any
 element of $F_2$. By  Lemma~\ref{lem:avila2} (iii), we have that
 $$
  |\Sigma(f_2)| < \E^{-\delta_1 \tilde{p} / 2}.
 $$
 Hence, \eqref{equa:opt} turns into
 $$
  \frac{1}{2} \delta_1 \ti{p} \geq \ti{p} \psi(\E^{-\delta_1 \tilde{p} / 2}).
 $$
 The claim now follows from the fact, that $\psi(x) \to 0$ as $x \to 0$.
\end{proof}

It now remains to construct the sequence of potentials $V_j$. We proceed by induction.
By possibly modifying $(\Omega,T)$, we can assume that
$$
 V_0(n) = f_0(T^n\omega)
$$
for some $f_0 \in C(\Omega)$. Assume now that, we are
given $V_1 = f_1 \circ T^n, \dots, V_{j-1} = f_{j-1} \circ T^n$ and we wish to construct
$V_{j} = f_j \circ T^n$. We can choose now
$\eps$ in the previous proposition to be the right hand side of \eqref{eq:Vjconv}.
We choose $\psi(x) = p_{j-1} \varphi(2 x)$ and the claim follows.

\end{document}